\documentclass[11pt,reqno]{amsart}

\usepackage[hmargin=1.2in]{geometry}
\usepackage{amsmath}
\usepackage{amsthm}
\usepackage{amscd}
\usepackage{xcolor}
\usepackage{amssymb}
\usepackage[pagebackref,hypertexnames=true]{hyperref}
\usepackage{stmaryrd}
\usepackage[all]{xy}

\theoremstyle{plain}
\newtheorem{thm}{Theorem}[section]
\newtheorem{cor}[thm]{Corollary}
\newtheorem{lemma}[thm]{Lemma}
\newtheorem{prop}[thm]{Proposition}

\theoremstyle{definition}

\newtheorem{defn}[thm]{Definition}

\newtheorem{remark}[thm]{Remark}

\newcommand{\f}{\varphi}
\newcommand{\ses}[3]{0 \to #1 \to #2 \to #3 \to 0}
\newcommand{\ilim}{\mathop{\varprojlim}}
\newcommand{\dlim}{\mathop{\varinjlim}}
\newcommand{\ilimone}{\mathop{\varprojlim}\,\hspace{-0.06in}^{1}}
\newcommand{\class}[1]{\llbracket #1 \rrbracket}
\newcommand{\hombbk}[2]{\text{\textnormal{Hom}}(\mathbb{K}_0(#1), \mathbb{K}_0(#2))}
\newcommand{\homk}[2]{\text{Hom}(K_0(#1), K_0(#2))}
\def\fa{\quad \forall}
\definecolor{MyGray}{rgb}{0.96,0.97,0.98}


\hypersetup{
	pdfauthor={Prahlad Vaidyanathan},
	colorlinks,
	linkcolor=blue,          
	citecolor=blue,        
	filecolor=magenta,      
	urlcolor=cyan,           
	pdfpagemode=UseNone
}

\title{E-theory for $C[0,1]$-algebras with finitely many singular points}
\author{Marius Dadarlat and Prahlad Vaidyanathan}

\begin{document}
\parindent 0pt

\thanks{M.D. was partially supported by NSF grant \#DMS--1101305.}
\begin{abstract}
We study the  E-theory group $E_{[0,1]}(A,B)$ for a class of C*-algebras over the unit interval with finitely many singular points, called elementary $C[0,1]$-algebras. We use results on  E-theory over non-Hausdorff spaces to describe $E_{[0,1]}(A,B)$ where $A$ is a sky-scraper algebra. Then we compute $E_{[0,1]}(A,B)$ for two elementary $C[0,1]$-algebras in the case where the fibers $A(x)$ and $B(y)$ of $A$ and $B$  are such that $E^1(A(x),B(y)) = 0$ for all $x,y\in [0,1]$.  This result  applies whenever the fibers satisfy the UCT, their $K_0$-groups are free of finite rank and their $K_1$-groups are zero. In that case we show that $E_{[0,1]}(A,B)$ is isomorphic to $\hombbk{A}{B}$, the  group of morphisms of the K-theory sheaves of $A$ and $B$. As an application, we give a streamlined partially new  proof of a classification result due to the first author and Elliott.
\end{abstract}
\maketitle

\section{Introduction}

 A deep isomorphism theorem of Kirchberg \cite{kirchberg} states that if $A$ and $B$ are strongly purely infinite, stable, nuclear, separable C*-algebras with primitive ideal spectrum homeomorphic to the same space $X$, then $A \cong B$ if and only the group $KK_X(A, B)$ contains an invertible element, where $KK_X(A,B)$ denotes the bivariant K-theory for C*-algebras over a space $X$. This leads to the question of computing $KK_X(A, B)$ and finding simpler invariants to understand this object. In the present paper we focus mainly on the case when $X$ is the unit interval. Recall that a C*-algebra over a locally compact, Hausdorff space $X$ is one that carries an essential central  action of $C_0(X)$, and by the Dauns-Hofmann theorem, every  C*-algebra with a Hausdorff spectrum $X$ can be thought of as a continuous $C_0(X)$-algebra (see Definition \ref{defn: cx_algebra}).\\

In this paper, we obtain information on $KK_X(A,B)$ using the $E$-theory groups $E_X(A,B)$, \cite{park_trout},\cite{mdd_meyer}.
 It is known \cite[Theorem 4.7]{park_trout} that $E_X(A,B)$ coincides with $KK_X(A,B)$ when $X$ is a locally compact Hausdorff space and $A$ is a separable continuous nuclear $C_0(X)$-algebra. Furthermore, the fact that E-theory satisfies excision for all extensions of $C_0(X)$-algebras enables us to compute the $E_{[0,1]}$-group for a class of elementary $C[0,1]$-algebras using the E-theory of their fibers and the E-theory classes of the connecting maps.
  We make crucial use of the generalization of $E$-theory to C*-algebras over non-Hausdorff spaces, as developed in  \cite{mdd_meyer}. Elementary $C[0,1]$-algebras, studied in \cite{mdd_fiberwise}, \cite{mdd_elliott} and \cite{mdd_pasnicu}, act as basic building blocks for more complex $C[0,1]$-algebras (see \cite[Theorem 6.2]{mdd_elliott}) and this paper should be viewed as a stepping stone towards that more general situation.
 As a first application of our calculations, we give a streamlined proof of the main result of \cite{mdd_elliott}, see Theorem~\ref{thm:DE}.
\\

One reason which makes the computation of the  $KK_X$-groups difficult is the prevalence of non-semisplit extensions over $X$. To illustrate this point, let us mention that the exact sequence of $C[0,1]$-algebras
 $0\to C_0[0,1)\to C[0,1]\to \mathbb{C}\to 0$
is not semisplit  over $X=[0,1]$. This is more than a technical nuisance, since
as pointed out in \cite[Remarques 1]{bauval},   a six-term  sequence of the form
\[
  \xymatrix{
    KK_X^0(\mathbb{C},D)\ar[r]&
    KK_X^0(C[0,1],D)\ar[r]&
    KK_X^0(C_0[0,1),D)\ar[d]\\
    KK_X^1(C_0[0,1),D)\ar[u]&
    KK_X^1(C[0,1],D)\ar[l]&
    KK_X^1(\mathbb{C},D)\ar[l]
  }
  \]
cannot be exact  for $D=C_0[0,1)$. Indeed, after computing each term one gets:
\[
  \xymatrix{
    0\ar[r]&
    0\ar[r]&
    \mathbb{Z}\ar[d]\\
    0\ar[u]&
    0\ar[l]&
    0\ar[l] {}.
  }
  \]
In contrast, the corresponding six-term exact sequence in $E^*_X$ is exact.  This property plays an important role
in our investigation of $C[0,1]$-algebras with finitely many singular points.

The paper is organized as follows: In Section \ref{sec: e_theory}, we revisit the construction of  E-theory for C*-algebras over a space $X$ and establish some preliminary lemmas. In Section \ref{sec: skyscraper}, we use the results in \cite{mdd_meyer} to describe $E_X(A,B)$ where $A$ is a sky-scraper algebra (Theorem \ref{thm: skyscraper}).
Section \ref{sec: elementary} contains the main result of this paper (Theorem \ref{thm: ex}), where we compute $E_X(A,B)$  for two elementary $C[0,1]$-algebras $A$ and $B$ whose fibers  satisfy the condition $E^1(A(x),B(y)) = 0$ for all $x,y\in [0,1]$.  In Section \ref{sec: k_theory}, we apply these results to the case where the fibers satisfy the UCT of \cite{rosenberg-schochet}, whose $K_0$-group is free of finite rank, and whose $K_1$-group is zero. In this case, we show (Theorem \ref{thm: uct}) that the natural map $\Gamma_{A,B}: E_X(A,B) \to \hombbk{A}{B}$ is an isomorphism. The fact that this map is surjective was proved, through different means, in \cite{mdd_elliott}. The merit of our arguments is that, not only were we able to show that the map is also injective, but we showed it using more natural methods. Also, the fact that we do  not require the UCT in Theorem \ref{thm: ex} is a marked difference from \cite{mdd_elliott}. \\

The authors are thankful to the referee for the suggestion to revisit the classification result of \cite{mdd_elliott}.

\section{E-theory over a topological space} \label{sec: e_theory}

We recall some basic definitions and facts from \cite{mdd_meyer}. \\

Let $A$ be a C*-algebra, and let Prim($A$) denote its primitive ideal space equipped with the hull-kernel topology. Let $\mathbb{I}(A)$ be the set of ideals of $A$ partially ordered by inclusion. For a space $X$, let $\mathbb{O}(X)$ be the set of open subsets of $X$ partially ordered by inclusion. Then there is a canonical lattice isomorphism \cite[$\mathsection$ 3.2]{dixmier}
\begin{equation}\label{eqn:lattice_iso}
\mathbb{O}(\text{Prim}(A)) \cong \mathbb{I}(A), \quad U \mapsto \bigcap \{\mathfrak{p}: \mathfrak{p} \in \text{Prim}(A)\setminus U \}
\end{equation}
\begin{defn} \label{defn: cx_algebra}
Let $X$ be a second countable topological space. A \emph{C*-algebra over $X$} is a C*-algebra $A$ together with a continuous map $\psi: \text{Prim}(A) \to X$. If in addition $X$ is locally compact and Hausdorff, then this is equivalent to a *-homomorphism from $C_0(X)$ to the center of the multiplier algebra of $A$ such that $C_0(X)A=A$. In this case, $A$ is called a $C_0(X)$-algebra. \\

For $U \subset X$ open, let $A(U) \in \mathbb{I}(A)$ be the ideal that corresponds to $\psi^{-1}(U) \in \mathbb{O}(\text{Prim}(A))$ under the isomorphism \eqref{eqn:lattice_iso}.  For $F \subset X$ closed, let $A(F) = A/A(X\setminus F)$. Both $A(U)$ and $A(F)$ are C*-algebras over $X$. \\

If $X$ is Hausdorff, $A(U) = C_0(U)A$. We write $\pi_x$ for the quotient map $A \to A(\{x\})$ and we say that $A$ is a continuous $C_0(X)$-algebra if the function $x\mapsto \|\pi_x(a)\|$ is continuous for all $a\in A$.
\end{defn}

\begin{defn}
An asymptotic morphism between two C*-algebras $A$ and $B$ is a family of maps $\f_t : A\to B$, for $t \in T := [0,\infty)$ such that $t \mapsto \f_t(a)$ is a bounded continuous function from $T$ to $B$ for each $a \in A$, and
$$
\f_t(a^{\ast} + \lambda b) - \f_t(a)^{\ast} - \lambda\f_t(b) \quad \text{ and }\quad \f_t(ab) - \f_t(a)\f_t(b)
$$
converge to $0$ in the norm topology as $t \to \infty$ for each $a, b \in A$ and $\lambda \in \mathbb{C}$. \\

Equivalently, an asymptotic morphism can be viewed as a map $\f: A \to C_b(T,B)$ that induces a *-homomorphism
$$
\hat{\f} : A \to B_{\infty} = C_b(T,B)/C_0(T,B).
$$

Two asymptotic morphisms $\f_t$ and $\f^{\prime}_t$ are called equivalent (in symbols, $\f_t \sim \f^{\prime}_t$) iff $\hat{\f} = \hat{\f^{\prime}}$, ie. $\f_t(a) - \f^{\prime}_t(a) \to 0$ as $t \to \infty$ for each $a \in A$. \\
\end{defn}

\begin{defn}
Let $A$ and $B$ be C*-algebras over a second countable topological space $X$. A *-homomorphism $\theta : A\to B$ is called $X$-equivariant if $\theta$ maps $A(U)$ into $B(U)$ for all open sets $U \subset X$. \\

An asymptotic morphism $\f : A\to C_b(T,B)$ is called approximately $X$-equivariant if, for any open set $U \subset X$,
\begin{equation}\label{eq:equiv}
\f(A(U)) \subset C_b(T,B(U)) + C_0(T,B).
\end{equation}
If $X$ is  second countable, then an asymptotic morphism $\f : A\to C_b(T,B)$ is  approximately $X$-equivariant if and only if it satisfies \eqref{eq:equiv} for only those open sets $U_n$ in a countable subbasis $(U_n)_n$ of the topology of $X$.

If in addition $X$ is a locally compact Hausdorff space, then by \cite[Lemma 2.11]{mdd_meyer} $\f$ is an approximately $X$-equivariant morphism if and only if $\f$ is an asymptotic $C_0(X)$-morphism in the sense of \cite[Definition 3.1]{park_trout}. i.e. $\f(fa)-f\f(a)\in C_0(T,B)$ for all $f\in C_0(X)$ and $a\in A$.
\end{defn}

\begin{lemma}\label{lem: restriction}
Let $A$ and $B$ be C*-algebras over a second countable topological space $X$ and let $\psi_t:A\to B$ be an approximately $X$-equivariant asymptotic morphism. For any open set $U \subset X$, there is an approximately $X$-equivariant asymptotic morphism $\psi_t^U : A \to B$ such that $\psi_t^U \sim \psi_t$ and
$$
\psi_t^{U}(A(U)) \subset B(U) \fa t \in T.
$$
\end{lemma}
\begin{proof} For an open subset $U$ of $X$ let

$$
B_{\infty}(U):=\frac{C_b(T,B(U))+C_0(T,B)}{C_0(T,B)}.
$$
Note that $B_{\infty}(X)=B_{\infty}$. Let $\pi:C_b(T,B)\to B_{\infty}$ be the quotient map.
 Let $s_U: B_\infty(U)\to C_b(T,B(U))$ be a set theoretic section of the surjective map $\pi_U :C_b(T,B(U))\to B_\infty(U)$ obtained by restricting $\pi$ to $C_b(T,B(U))$.
Extend $s_U$ to a section $s:B_\infty \to C_b(T,B)$ of $\pi$. Then $\psi^{U}:=s\circ \hat \psi$ is an asymptotically $X$-equivariant asymptotic morphism equivalent to $\psi$ since $\hat \psi^{U}=\hat \psi$.  Moreover, we have that
$$\psi^{U}(A(U))=s(\hat \psi(A(U)))\subset C_b(T,B(U))$$
since $\hat\psi(A(U))\subset B_{\infty}(U)$ as a consequence of the assumption that $\psi(A(U)) \subset C_b(T,B(U)) + C_0(T,B)$.

\end{proof}

\begin{defn}
A homotopy of asymptotic morphisms from $A$ to $B$ is an asymptotic morphism from $A$ to $C([0,1],B)$. Let $\llbracket A,B\rrbracket_X$ denote the set
of homotopy classes of approximately $X$-equivariant asymptotic morphisms from A to B, and let $\class{\psi_t}$ denote the homotopy class of an approximately $X$-equivariant asymptotic morphism $\psi_t:A \to B$. \\

It is immediate that equivalent asymptotic morphisms are homotopic.
\end{defn}

\begin{defn}
Let $X$ be a second-countable topological space and let $A$ and $B$ be separable C*-algebras over $X$. Define
\begin{equation*}
\begin{split}
E_X(A,B) = E_X^0(A,B) &= \llbracket SA\otimes \mathcal{K}, SB\otimes \mathcal{K} \rrbracket_X \\
E_X^1(A,B) &= E_X(A,SB)
\end{split}
\end{equation*}
where $S$ denotes the suspension functor $SA = C_0(\mathbb{R},A)$. \\

By \cite[Theorem 2.25]{mdd_meyer}, there is a composition product
$$
E_X(A,B)\times E_X(B,C) \to E_X(A,C)
$$
and $E_X(\cdot, \cdot)$ is the universal half-exact, C*-stable homotopy functor on the category of separable C*-algebras over $X$. There are six-term exact sequences in each variable of $E^{\ast}_X(A,B)$. \\

Furthermore, if $X$ is a locally compact, Hausdorff space, then this definition of $E^{\ast}_X(A,B)$ coincides with that of Park and Trout (see \cite[Proposition 2.29]{mdd_meyer})
\end{defn}

With these definitions in place, Proposition \ref{prop: open_set} is now a simple consequence of Lemma \ref{lem: restriction}.

\begin{prop} \label{prop: open_set}
Let $X$ be a second-countable topological space and let $A$ and $B$ be separable C*-algebras over $X$. For any open set $U \subset X$, the inclusion map $i: B(U) \to B$ induces a natural isomorphism
$$
i_{\ast} : E_X^{\ast}(A(U), B(U)) \xrightarrow{\cong} E_X^{\ast}(A(U), B).
$$
\end{prop}

\begin{prop} \label{prop: restriction}
Let $X$ be a second-countable topological space and let $A$ and $B$ be separable C*-algebras over $X$.
If $U \subset Y \subset X$ with $U$ open and $Y$ has the induced topology, then
$$
E_X^{\ast}(A(U), B(U)) \cong E_Y^{\ast}(A(U), B(U)).
$$
\end{prop}
\begin{proof} Let us first observe that $A(U)$ and $B(U)$ are C*-algebras over $Y$. Indeed if $W$ is open in $Y$,
then $W=Y\cap V$ for some open subset $V$ of $X$ and $A(U)(W):=A(U\cap V)$.  It is then
clear that an asymptotic morphism $\f:A(U)\to C_b(T,B(U))$ is approximately $X$-equivariant if and only if it is approximately $Y$-equivariant. This concludes the proof.
\end{proof}


\section{Sky-Scraper Algebras} \label{sec: skyscraper}

In this section, we consider C*-algebras over a space $X$ with exactly one non-trivial fibre, called sky-scraper algebras. We use \cite[Theorem 3.2]{mdd_meyer} to exhibit a short-exact sequence that computes $E_X(A,B)$, where $A$ is a sky-scraper algebra. In the following section, we use this exact sequence to isolate those points  where a $C[0,1]$-algebra is not locally trivial.

\begin{defn}[Sky-scraper Algebra] \label{ex: skyscraper}
Let $D$ a separable C*-algebra and let $X$ be a topological space. Fix a point $x\in X$ and define $i_x(D)$ to be the C*-algebra $D$ regarded as a C*-algebra over $X$ by setting
$$
i_x(D)(U) = \begin{cases}
			D & : x \in U \\
			0 & : x \notin U.
			\end{cases}
$$
$i_x(D)$ is called the sky-scraper algebra with fiber $D$ at the point $x$.
If $X$ is locally compact, the corresponding action  of $C_0(X)$ is given by
$\iota(f)(d) := f(x)d$.
\end{defn}

\begin{thm} \label{thm: skyscraper}
Let $i_x(D)$ be a sky-scraper algebra as in Definition \ref{ex: skyscraper} and assume that $X$ is second countable.
 Let $\{U_n\}_n$ be a neighbourhood basis of open neighbourhoods of $x$ such that
$U_{n+1} \subset U_n \fa n \in \mathbb{N}.$
Then, for any separable C*-algebra $B$ over $X$, there is a short exact sequence
$$
\ses{\ilimone{E^{\ast + 1}(D, B(U_n))}}{E_X^{\ast}(i_x(D), B)}{\ilim{E^{\ast}(D, B(U_n))}}
$$
\end{thm}
\begin{proof}
Let $Y := (X, \tau)$ be a topological space whose underlying space is $X$, but whose topology $\tau$ is the topology generated by the sets $\{U_n\}_n$. We claim that
\begin{equation}\label{eq:finite}
E_X^{\ast}(i_x(D), B) = E_Y^{\ast}(i_x(D), B)
\end{equation}
To see this, consider an asymptotic morphism $\psi_t : i_x(D) \to B$ and the induced map $\psi : i_x(D) \to C_b(T,B)$ (we do not use suspensions for ease of notation, but it is clear that the same argument holds with suspensions). From the definition of $i_x(D)$,
we see that $\psi_t$ is approximately $X$-equivariant if and only if
$$
\psi(D) \subseteq C_b(T,B(U)) + C_0(T, B) \fa \,\, U \text{open}\, \subset X, \, x \in U.
$$
On the other hand $\psi_t$ is approximately $Y$-equivariant, if and only if
$$
\psi(D) \subseteq C_b(T,B(U_n)) + C_0(T, B) \fa \, n \in \mathbb{N}
$$
But if $U \subset X$ is open and $x \in U$, then there is $n_0 \in \mathbb{N}$ such that $U_{n_0} \subset U$. Thus, if $\psi_t$ is asymptotically $Y$-equivariant, then
$$
\psi(D) \subseteq C_b(T,B(U_{n_0})) + C_0(T, B) \subseteq C_b(T,B(U)) + C_0(T, B)
$$
and hence $\psi_t$ will be approximately $X$-equivariant as well. The converse is obvious since every open set in $Y$ is already open in $X$.
Since the same argument applies to homotopies $\Phi_t : i_x(D) \to C([0,1],B)$, we obtain \eqref{eq:finite}.
With a view to apply the $\ilimone$-sequence from \cite[Theorem 3.2]{mdd_meyer}, we define
$
X_n := (X, \tau_n),$  where the topology $\tau_n$ is generated by the sets $\{U_1, U_2, \ldots, U_n\}$.
 We now claim that
$$
E_{X_n}^{\ast}(i_x(D), B) \cong E^{\ast}(D, B(U_n)).
$$
As before (omitting suspensions), we consider an asymptotic morphism $\psi_t : i_x(D) \to B$, and note that, since $U_n \subset U_i$ for each $i \leq n$, $\psi_t$ is  approximately  $X_n$-equivariant if and only if
$$
\psi(D) \subseteq C_b(T,B(U_n)) + C_0(T, B).
$$
Since $U_n \subset X_n$ is open, we may apply Lemma \ref{lem: restriction} to obtain a map
$$
\eta : E_{X_n}(i_x(D), B) \to E(D, B(U_n)), \qquad \class{\psi_t} \mapsto \class{\psi_t^{U_n}}.
$$
We claim that $\eta$ is bijective: Suppose that $\class{\psi_t^{U_n}} = \class{\f_t^{U_n}}$ in $E(D, B(U_n))$. Then there is a homotopy $\Phi_t: D \to C[0,1]\otimes B(U_n)$ connecting $\psi_t^{U_n}$ to $\f_t^{U_n}.$ Since
$$
\Phi_t(D) \subset C[0,1]\otimes B(U_n) \subset C[0,1]\otimes B,
$$
$\Phi_t$ is an asymptotically $X_n$-equivariant map from $i_x(D)$ to $C[0,1]\otimes B$ connecting $\psi_t^{U_n}$ and $\f_t^{U_n}$. But by Lemma \ref{lem: restriction}, $\psi_t \sim \psi_t^{U_n}$ and $\f_t \sim \f_t^{U_n}$, and hence $\class{\psi_t} = \class{\f_t}$ in $E_{X_n}(i_x(D), B)$ and hence $\eta$ is injective. \\
For surjectivity, we observe that any given asymptotic morphism $\f_t : D \to B(U_n)$,
can be viewed as an  $X_n$-equivariant asymptotic morphism $\psi_t:i_x(D)\to B$.
We are now in a position to complete the proof. Since the collection $\{U_n\}$ forms a countable basis for the topological space $Y,$ we may apply \cite[Theorem 3.2]{mdd_meyer} to obtain a short exact sequence
$$
\ses{\ilimone{E^{\ast + 1}_{X_n}(i_x(D),B)}}{E_Y^{\ast}(i_x(D),B)}{\ilim{E^{\ast}_{X_n}(i_x(D),B)}}.
$$
By our earlier identifications, this reduces to
$$
\ses{\ilimone{E^{\ast + 1}(D, B(U_n))}}{E_X^{\ast}(i_x(D), B)}{\ilim{E^{\ast}(D, B(U_n))}}. \qedhere
$$
\end{proof}

We now list some corollaries of Theorem \ref{thm: skyscraper} which will be useful to us in the next section.
\begin{cor} \label{cor: cor_sky_1}
 If $U \subset X$ is an open set such that $x \notin \overline{U}$, then
$$
E_X^{\ast}(i_x(D), B(U)) = 0.
$$
\end{cor}
\begin{proof}
There exists an open neighbourhood $O$ of $x$ such that $O\cap U = \emptyset$. Consider a sequence of neighbourhoods $(U_n)_n$ of $x$ as in Theorem \ref{thm: skyscraper} such that $U_1 = O$. Then
$$
E^{\ast}(D, B(U)(U_n)) = E^{\ast}(D, B(U\cap U_n)) \cong 0,
$$
and the result follows from Theorem \ref{thm: skyscraper}.
\end{proof}

\begin{cor} \label{cor: cor_sky_2}
Let $A$, $B$, $D$, $\{U_n\}_n$ and $x$ be as in Theorem~\ref{thm: skyscraper}.
Suppose that the all the inclusions $B(U_{n+1})\hookrightarrow B(U_{n})$ are equivalences in E-theory. Then, for any fixed $k \in \mathbb{N}$:
$$
E_X^{\ast}(i_x(D), B) \cong E^{\ast}(D, B(U_k)).
$$
\end{cor}

\begin{proof}

The assumption that the inclusion maps $B(U_{n+1}) \hookrightarrow B(U_n)$ is an equivalence in E-theory
implies that
 $\ilim{E^{\ast}(D, B(U_n))} \cong E^{\ast}(D, B(U_k))$ and
$
\ilimone{E^{\ast+1}(D, B(U_n))} \cong 0.
$
The conclusion follows now from Theorem \ref{thm: skyscraper}.
\end{proof}

\begin{cor} \label{cor: e_trivial}
Let $U \subset [0,1]$ be an open interval. For any two separable C*-algebras $D, H$
$$
E_{[0,1]}^{\ast}(C_0(U)\otimes D, C_0(U)\otimes H) \cong E^{\ast}(D,H)
$$
\end{cor}
\begin{proof}
Suppose that $U = (a,b)$, then by Proposition \ref{prop: restriction}, we may assume without loss of generality that $a = 0$ and $b = 1$. Now by Theorem \ref{prop: open_set}
$$
E_{[0,1]}^{\ast}(C_0(0,1)\otimes D, C_0(0,1)\otimes H) \cong E_{[0,1]}^{\ast}(C_0(0,1)\otimes D, C[0,1]\otimes H)
$$
Now consider the short exact sequence
$$
\ses{C_0(0,1)\otimes D}{C[0,1]\otimes D}{i_0(D)\oplus i_1(D)}.
$$
By Corollary \ref{cor: cor_sky_1},
$$
E_{[0,1]}^{\ast}(i_0(D), C[0,1]\otimes H) \cong E^{\ast}(D, C_0[0,1)\otimes H) \cong 0
$$
since $C_0[0,1)\otimes H$ is contractible. Similarly
$$
E_{[0,1]}^{\ast}(i_1(D), C[0,1]\otimes H) \cong 0.
$$
Hence by using the six-term exact sequence in the first variable of $E^*_{[0,1]}$ and \cite[Lemma~2.30]{mdd_meyer} we obtain
\begin{equation*}
\begin{split}
E_{[0,1]}^{\ast}(C_0(0,1)\otimes D, C_0(0,1)\otimes H) & \cong E_{[0,1]}^{\ast}(C[0,1]\otimes D, C[0,1]\otimes H) \\
													   & \cong E^{\ast}(D, C[0,1]\otimes H) \\
													   & \cong E^{\ast}(D, H)
\end{split}
\end{equation*}
Now suppose $U = [0,a)$ or $U = (b,1]$, and since the proofs are identical, we assume that $U = [0,a)$. Furthermore, by using Proposition \ref{prop: restriction}, we may assume without loss of generality that $a = 1$. Now the proof is identical to the first part, except that we use the short exact sequence
$$
\ses{C_0[0,1)\otimes D}{C[0,1]\otimes D}{i_1(D)}
$$
instead.
\end{proof}


\section{Elementary $C[0,1]$-algebras} \label{sec: elementary}

A $C[0,1]$-algebra is said to be locally trivial at a point $x \in [0,1]$ if there is an open neighborhood $U$ of $x$, and a C* algebra $D$ such that $A(U) \cong C_0(U)\otimes D$. If $A$ is not locally trivial at $x$, we say that $x$ is a singular point of $A$. \\

By an \emph{elementary $C[0,1]$-algebra}, we mean an algebra which is locally trivial at all but finitely many points and moreover the algebra has a specific structure at the singular points as described below in Definition~\ref{def:elementary-1}. The importance of such algebras comes from the following theorem due to the first author and Elliott.

\begin{thm}\label{thm:structure}\cite[Theorem 6.2]{mdd_elliott}
Let $\mathcal{C}$ be a class of unital semi-projective Kirchberg algebras. Let $A$ be a separable unital continuous $C[0,1]$-algebra such that all of its fibers are inductive limits of sequences in $\mathcal{C}$. Then, there exists an inductive system $(A_k, \varphi_k)$ consisting of unital elementary $C[0,1]$-algebras with fibers in $\mathcal{C}$ and unital morphisms of $C[0,1]$-algebras $\varphi_k \in \text{\textnormal{Hom}}(A_k,A_{k+1})$ such that
$$
A \cong \dlim (A_k, \varphi_k)
$$
A similar result is valid if one assumes that all the C*-algebras in $\mathcal{C}$ are stable rather than unital.
\end{thm}
Theorem~\ref{thm:structure} applies to all continuous $C[0,1]$-algebras whose fibers are Kirchberg algebras
satisfying the UCT and having torsion free $K_1$-groups.
Furthermore, by \cite[Theorem 2.5]{mdd_fiberwise}, any separable nuclear continuous $C[0,1]$-algebra over $[0,1]$ is $KK_{[0,1]}$-equivalent to a separable continuous unital $C[0,1]$-algebra whose fibers are Kirchberg algebras. Thus, the elementary $C[0,1]$-algebras are basic building blocks (in a K-theoretical sense) of all continuous $C[0,1]$-algebras. \\

Elementary $C[0,1]$-algebras are given by the following data: Let $\mathcal{F}$ be a fixed class of separable C*-algebras. Let $X = [0,1]$, and consider a partition $\mathcal{P}$ of $[0,1]$ given by
$$
0 = x_0 <x_1 < \ldots < x_n < x_{n+1} = 1
$$
Write $F := \{x_1, x_2, \ldots, x_n\}$. We define a $C[0,1]$-algebra which is locally trivial at all points except possibly this finite set $F$ and has fibers in the class $\mathcal{F}$. \\

Suppose that we are given C*-algebras
$$
\{D_1, D_2, \ldots, D_n, H_1, H_2, \ldots, H_{n+1}\} \subset \mathcal{F}
$$
and *-homomorphisms
$$
\gamma_{i,0} : D_i \to H_i \qquad \gamma_{i,1} : D_i \to H_{i+1}
$$

Define
\begin{equation*}
\begin{split}
A = \begin{Bmatrix} ((d_1, d_2, \ldots, d_n), (h_1, h_2,\ldots, h_{n+1})) :& d_i \in D_i, h_i \in C[x_{i-1},x_i]\otimes H_i \\
& \text{ s.t. } \gamma_{i,j}(d_i) = h_{i+j}(x_i) \\
&\forall i \in \{1,2,\ldots, n\}, j\in \{0,1\} \end{Bmatrix}
\end{split}
\end{equation*}

In other words, $A$ is the pull-back of the diagram
$$
\begin{CD}
A @>>> \bigoplus_{i=1}^{n+1} C[x_{i-1},x_i]\otimes H_i \\
@VVV @VV{eval}V \\
\bigoplus_{i=1}^n D_i @> (\gamma_{i,0}, \gamma_{i,1}) >> \bigoplus_{i=1}^n H_i\oplus H_{i+1}
\end{CD}
$$
\begin{defn}\label{def:elementary-1}
A $C[0,1]$-algebra $A$ as above that is associated to the partition $\mathcal{P}$ with fibers in $\mathcal{F}$ and which satisfies the condition that
 \emph{for each $i \in \{1,2,\ldots, n\}$, either $\gamma_{i,0}$ or $\gamma_{i,1}$ is the identity map,}
 is called an elementary $C[0,1]$-algebra.
 We denote the class of all such algebras by $\mathcal{E}(\mathcal{P}, \mathcal{F})$.
Note that if $\mathcal{P}_1$ is a refinement of $\mathcal{P}_2$, then $\mathcal{E}(\mathcal{P}_2, \mathcal{F}) \subset \mathcal{E}(\mathcal{P}_1, \mathcal{F})$ since we may add singularities by choosing the maps $\gamma_{i,j}$ to be the identity maps. We define
$$
\mathcal{E}(\mathcal{F}) := \bigcup_{\mathcal{P}}\mathcal{E}(\mathcal{P}, \mathcal{F})=\text{class of all  elementary $C[0,1]$-algebra with fibers in } \mathcal{F}.
$$
When we write $A, B \in \mathcal{E}(\mathcal{F})$, we implicitly mean that we are choosing a common partition $\mathcal{P}$ as above.
\end{defn}


We are now concerned with computing $E_X(A,B)$ for $A, B \in \mathcal{E}(\mathcal{F})$. In order to simplify our future work, we fix $A$ as above, and define $B$ as
\begin{equation*}
\begin{split}
B = \begin{Bmatrix} ((d_1', d_2', \ldots, d_n'), (h_1', h_2',\ldots, h_{n+1}')) :& d_i' \in D_i', h_i' \in C[x_{i-1},x_i]\otimes H_i' \\
& \text{ s.t. } \gamma_{i,j}'(d_i') = h_{i+j}'(x_i) \\
& \forall i\in \{1,2,\ldots, n\}, j\in \{0,1\} \end{Bmatrix}
\end{split}
\end{equation*}

In other words, the fibers of $A$ will be $D_i$ or $H_j$ and the connecting maps will be $\gamma_{i,j}$, while the corresponding fibers of $B$ will be $D_i'$ or $H_j'$, and the connecting maps of $B$ will be $\gamma_{i,j}'$. \\

Now consider the partition $\mathcal{P}$ as above, and write
$$
U = [0,1]\setminus F = \bigsqcup_{i=1}^{n+1} U_i
$$
where $U_1 = [0,x_1), U_{n+1} = (x_n,1]$ and $U_i = (x_{i-1}, x_i)$ for $2 \leq i \leq n$. The short exact sequence
\begin{equation} \label{eqn: ses_a}
\ses{A(U)}{A}{A(F)}
\end{equation}
yields a long exact sequence in E-theory

\begin{equation} \label{eqn: les_eab}
\begin{CD}
E_X(A(F),B) @>>> E_X(A, B) @>>> E_X(A(U), B) \\
@AAA				@.				@VV\delta V \\
E_X^1(A(U), B) @<<< E_X^1(A,B) @<<< E_X^1(A(F), B)
\end{CD}
\end{equation}

whose boundary map we denote by
$$
\delta : E_X(A(U),B) \to E_X^1(A(F),B)
$$
As we will show later, this map $\delta$ holds the key to understanding $E_X(A,B)$. We begin by identifying the domain of this map.

\begin{lemma} \label{lem: au_b}
The inclusion map $B(U) \hookrightarrow B$ induces an isomorphism
$$
E_X(A(U),B)\cong E_X(A(U),B(U)) \cong \bigoplus_{i=1}^n E(H_i,H_i')
$$
\end{lemma}
\begin{proof}
The first isomorphism follows from Proposition~\ref{prop: open_set}. Furthermore, since $E_{X}(A(U_i),B(U_j))=0$ if $i\neq j$, by additivity of $E_X$ in each variable:
\begin{equation*}
\begin{split}
E_X(A(U),B(U)) &\cong \bigoplus_{i=1}^{n+1} E_X(A(U_i), B(U_i)) \\
&\cong E_{\overline{U_i}}(A(U_i),B(U_i) \qquad\qquad\text{(by Proposition \ref{prop: restriction})}\\
&\cong \bigoplus_{i=1}^{n+1} E(H_i,H_i') \qquad\qquad\quad\text{(by Corollary \ref{cor: e_trivial})}
\end{split}
\end{equation*} \qedhere
\end{proof}

\begin{remark} \label{rem: delta_u}
Consider the inclusion $B(U)\hookrightarrow B$ and the induced commutative diagram
$$
\begin{CD}
E_X(A(U),B)) @> \delta >> E_X^1(A(F),B) \\
@AAA @AA \iota A \\
E_X(A(U),B(U)) @> \Delta_A >> E_X^1(A(F),B(U))
\end{CD}
$$
By Lemma \ref{lem: au_b}, the vertical map on the left is an isomorphism, so
$$
\ker(\delta) \cong \ker(\iota\circ \Delta_A).
$$
We now compute the map $\Delta_A$. Consider the short exact sequence
$$
\ses{A(U)}{A}{A(F)}
$$
and the boundary element in E-theory obtained from this sequence,
$$
\delta_A \in E_X^1(A(F),A(U))
$$
and note that $\Delta_A$ is given by multiplication by this element
$$
E_X^1(A(F),A(U)) \ni \delta_A \times E_X(A(U),B(U)) \xrightarrow{\Delta_A} E_X^1(A(F),B(U))
$$
The next two lemmas help us compute this map.
\end{remark}

\begin{lemma} \label{lem: af_bu}
$
E_X^1(A(F),B(U)) \cong \bigoplus_{i=1}^n (E(D_i,H_i')\oplus E(D_i,H_{i+1}')).
$
\end{lemma}
\begin{proof}
By additivity of $E_X$
\begin{equation*}
\begin{split}
E_X^1(A(F),B(U)) &\cong \bigoplus_{i=1}^n E_X^1(A(x_i), B(U)) \\
&\cong \bigoplus_{i=1}^n E_X^1(i_{x_i}(D_i), B(U_i\cup U_{i+1})) \qquad\text{(by Corollary \ref{cor: cor_sky_1})} \\
&\cong \bigoplus_{i=1}^n E^1(D_i,B(U_i\cup U_{i+1})) \qquad\qquad \text{(by Corollary \ref{cor: cor_sky_2})} \\
&\cong \bigoplus_{i=1}^n (E^1(D_i,C_0(U_i)\otimes H_i')\oplus E^1(D_i,C_0(U_{i+1})\otimes H_{i+1}') \\
&\cong \bigoplus_{i=1}^n (E(D_i,H_i')\oplus E(D_i,H_{i+1}')).
\end{split}
\end{equation*} \qedhere
\end{proof}

\begin{lemma} \label{lem: delta_a}
$
E_X^1(A(F),A(U)) \cong \bigoplus_{i=1}^n \left(E(D_i,H_i) \oplus E(D_i,H_{i+1})\right)
$
and under this isomorphism
\begin{equation} \label{eqn: delta_a}
\delta_A \mapsto (-\class{\gamma_{i,0}}, \class{\gamma_{i,1}})_{i=1}^n
\end{equation}
\end{lemma}
\begin{proof}
The isomorphism from the statement follows from Lemma \ref{lem: af_bu} applied for $B=A$.
In order to compute the image of $\delta_A$ under this isomorphism,
we need the following notation:
 $U_{1,0} = [0,x_1], U_{n+1,1} = [x_n,1]$ and $U_{i,0} = (x_{i-1}, x_i]$, $U_{i,1} = [x_i,x_{i+1} )$ for $2 \leq i \leq n$.
For each $i\in \{1,...,n\}$ and $j\in\{0,1\}$ consider the extension of $C[0,1]$-algebras
\begin{equation}\label{ds}
\ses{A(U_{i+j})}{A(U_{i,j})}{A(x_i)}
\end{equation}
and the corresponding element $\delta_{i,j}$ that belongs to $E_X^1(A(x_i), A(U_{i+j}))$ which we may regard as direct summand of $E^1_X(A(F),A(U))$.
We are going to show that
$$\delta_A=\bigoplus_{i=1}^n (\delta_{i,0}\oplus \delta_{i,1})\in \bigoplus_{i=1}^n \left(E_X^1(A(x_i), A(U_{i} ))\oplus E_X^1(A(x_i), A(U_{i+1} ))\right).$$
To this purpose we will write explicitly an expression for the Connes-Higson  asymptotic morphism $(\gamma_t)_{t\in [0,1)}$
that defines $\delta_A$, see \cite[Prop.~2.23]{mdd_meyer}. Let $(u_i^t)_{t\in [0,1)}$ be a contractive positive approximate unit of $C_0(U_i)$.
For each $i$, choose two continuous maps $\eta_{i,0}\in C_0(x_{i-1},x_i]$ and  $\eta_{i,1}\in C_0[x_i,x_{i+1})$ such that they assume values in $[0,1]$,
 they are equal to $1$ on a neighborhood of $x_i$  and such  that $\eta_{i,1}\eta_{i+1,0}=0$ for $1\leq i <n$.
 It follows that we have the following asymptotic expression for $(\gamma_t)_{t\in [0,1)}:C_0(0,1)\otimes A(F)\to A(U)$,
 \[\gamma_t (f \otimes (d_i)_{i=1}^n)=\sum_{i=1}^n \, \left( f(u_i^t)\eta_{i,0} \otimes \gamma_{i,0}(d_{i}) +    f(u_{i+1}^{t})\eta_{i,1} \otimes \gamma_{i,1}(d_{i})\right). \]

It is now clear that $\gamma_t$ decomposes in orthogonal sum of  components

$\gamma^{i,j}_t(f\otimes d_i)=f(u_{i+j}^t)\eta_{i,j} \otimes \gamma_{i,j}(d_i)$, $1\leq i \leq n$, $0\leq j \leq 1.$ But we now recognize $\gamma^{i,j}_t$ as representing the Connes-Higson asymptotic morphism
defined by the extension \eqref{ds}.  Next we are going to identify its class. We focus on the point $x_i$, and consider the map of extensions
$$
\begin{CD}
0 @>>> C_0(U_{i+j})\otimes H_{i+j} @>>> A({U_{i,j}}) @>>> i_{x_i}(D_i) @>>> 0 \\
@. @VV = V @VVV @VV \gamma_{i,j} V @. \\
0 @>>> C_0(U_{i+j})\otimes H_{i+j} @>>> C_0({U_{i,j}})\otimes H_{i+j} @>>> i_{x_i}(H_{i+j}) @>>> 0
\end{CD}
$$
We apply the functor $E_X(\cdot, A(U_{i+j}))$ to this sequence, and consider the relevant part of the resulting commutative diagram
$$
\begin{CD}
E_X(C_0(U_{i+j})\otimes H_{i+j}, A(U_{i+j})) @> \delta^{i,j}_A >> E_X^1(i_{x_1}(D_1), A(U_{i+j})) \\
@A = AA @AA \gamma_{i,j}^{\ast} A \\
E_X(C_0(U_{i+j})\otimes H_{i+j}, A(U_{i+j})) @> \delta^{i,j} >> E_X^1(i_{x_1}(H_{i+j}), A(U_{i+j}))
\end{CD}
$$
where the map $\delta^{i,j}$ is given by multiplication by the boundary element
$$
\class{\delta_t} \in E_X^1(i_{x_i}(H_{i+j}), C_0(U_{i+j})\otimes H_{i+j}) \cong E^1(H_{i+j}, C_0(U_{i+j})\otimes H_{i+j}))
$$
(by Corollary \ref{cor: cor_sky_2}).  If  $j=0$, $\class{\delta_t}$ corresponds  under this isomorphism to the boundary map of the extension
$$
\ses{C_0(0,1)\otimes H_i}{C_0(0, 1]\otimes H_i}{H_i}
$$
 which can be identified with the element $-\class{\text{id}} \in E(H_i, H_i)$. This accounts for the  sign of the term $\class{\gamma_{i,0}}$ in the expression \eqref{eqn: delta_a}. \\

Similarly, if $j=1$, $\class{\delta_t}$
corresponds to the boundary map of the extension of C*-algebras
$$
\ses{C_0(0, 1)\otimes H_{i+1}}{C_0[0, 1)\otimes H_{i+1}}{H_{i+1}}
$$
which can be identified with the element $\class{\text{id}} \in E(H_{i+1},H_{i+1})$. This accounts for the difference in sign.
\end{proof}

The next lemma now follows from Remark \ref{rem: delta_u} and Lemmas \ref{lem: af_bu}, \ref{lem: delta_a}
\begin{lemma} \label{lem: delta_u} There is a commutative diagram
\[
  \begin{xy}
   \xymatrix{
E_X(A(U),B(U)) \ar[r]^{\Delta_A }\ar[d]_{\cong}& E_X^1(A(F),B(U)) \ar[d]_{\cong}\\
\bigoplus_{i=1}^{n+1} E(H_i,H_i') \ar[r] &\bigoplus_{i=1}^n (E(D_i,H_i')\oplus E(D_i,H_{i+1}'))
   }
  \end{xy}
\]
Under this isomorphisms, $\Delta_A$ corresponds to the map
$$
(\beta_i)_{i=1}^{n+1} \mapsto (- \beta_i \circ \class{\gamma_{i,0}}, \beta_{i+1}\circ \class{\gamma_{i,1}})_{i=1}^n
$$
\end{lemma}

\begin{remark} \label{rem: iota}
We saw in Remark \ref{rem: delta_u} that
$$
\ker(\delta) \cong \ker(\iota\circ \Delta_A)
$$
where $\iota : E_X^1(A(F),B(U)) \to E_X^1(A(F),B)$ is induced by the inclusion $B(U) \hookrightarrow B$. In order to compute this kernel, consider the following long exact sequence coming from the extension $\ses{B(U)}{B}{B(F)}:$
$$
\begin{CD}
E_X^1(A(F),B(U)) @> \iota >> E_X^1(A(F),B)) @>>> E_X^1(A(F),B(F)) \\
@A \Delta_B AA @. @VVV \\
E_X(A(F),B(F)) @<<< E_X(A(F),B) @<<< E_X(A(F),B(U))
\end{CD}
$$
Thus,
$$
\ker(\iota) = \text{Im}(\Delta_B)
$$
where $\Delta_B$ is given by multiplication by the boundary element
$$
\delta_B \in E_X^1(B(F),B(U))
$$
As in Lemma \ref{lem: delta_a}, we have
$$
E_X^1(B(F),B(U)) \cong \bigoplus_{i=1}^n E(D_i',H_i') \oplus E(D_i',H_{i+1}')
$$
and under this isomorphism
\begin{equation} \label{eqn: delta_b}
\delta_B \mapsto (-\class{\gamma_{i,0}'}, \class{\gamma_{i,1}'})_{i=1}^n
\end{equation}
\end{remark}

\begin{lemma} \label{lem: Delta}
\[
  \begin{xy}
   \xymatrix{
E_X^1(A(F),B(U)) \ar[r]^-{\cong }& \bigoplus_{i=1}^n (E(D_i,H_i')\oplus E(D_i,H_{i+1}'))  \\
E_X(A(F),B(F)) \ar[r]^{\cong}\ar[u]^{\Delta_B} &\bigoplus_{i=1}^n E(D_i,D_i')\ar[u]
   }
  \end{xy}
\]
and under these isomorphisms
$$
\Delta_B((\alpha_i)_{i=1}^n) = (-\class{\gamma_{i,0}'}\circ \alpha_i, \class{\gamma_{i,1}'}\circ \alpha_i)
$$
\end{lemma}
\begin{proof}
The map $\Delta_B$ is induced by the product
$$E_X(A(F),B(F))\times \delta_B \in E_X^1(B(F),B(U))\to E_X^1(A(F),B(U).$$
We have already described all the terms that appear in this composition.
\end{proof}
\begin{thm} \label{thm: delta}
For $A,B \in \mathcal{E}(\mathcal{F})$ as above
\begin{equation}\label{eqn: ex}
\begin{split}
\ker(\delta) = \begin{Bmatrix} (\beta_i) \in \bigoplus_{i=1}^{n+1} E(H_i,H_i') : & \exists (\alpha_i) \in \bigoplus_{i=1}^n E(D_i,D_i') \text{ s.t.} \\
				& \beta_i\circ \class{\gamma_{i,0}} = \class{\gamma_{i,0}'}\circ \alpha_i \text{ and } \\
				& \beta_{i+1} \circ \class{\gamma_{i,1}} = \class{\gamma_{i,1}'} \circ \alpha_i \end{Bmatrix}
\end{split}
\end{equation}
\end{thm}
\begin{proof}
By Remarks \ref{rem: delta_u} and \ref{rem: iota}
\begin{equation*}
\begin{split}
\ker(\delta) & \cong \ker(\iota\circ\Delta_A) \\
& \cong \{\beta \in E_X(A(U),B(U)) : \Delta_A(\beta) \in \ker(\iota)=\text{Im}(\Delta_B) \} \\
&\cong \{\beta \in E_X(A(U),B(U)) : \exists \alpha \in E_X(A(F),B(F)) \text{ s.t. } \Delta_A(\beta) = \Delta_B(\alpha)\}
\end{split}
\end{equation*}
The expression in \eqref{eqn: ex} now follows from the description of $\Delta_A$ and $\Delta_B$ from Lemmas \ref{lem: delta_u} and \ref{lem: Delta}.
\end{proof}
\begin{defn}\label{defn:classF}
We now specify a type of class $\mathcal{F}$ for which we can explicitly compute $E_X(A,B)$ for any $A,B \in \mathcal{E}(\mathcal{F})$ using the machinery developed above. Let $\mathcal{F}$ be a class of separable C*-algebras such that $E^1(D,D') = 0$ for all $D,D'\in \mathcal{F}$.
\end{defn}

\begin{lemma} \label{lem: af_b}
If $A,B \in \mathcal{E}(\mathcal{F})$ with $\mathcal{F}$ as in Definition~\ref{defn:classF}, then
$
E_X(A(F),B) = 0.
$
\end{lemma}
\begin{proof}
By the additivity of $E_X(\cdot, B)$
$$
E_X(A(F),B) \cong \bigoplus_{i=1}^n E_X(A(x_i), B) \cong \bigoplus_{i=1}^n E_X(i_{x_i}(D_i), B).
$$
Choose $\epsilon > 0$ small enough so that $(x_i-\epsilon, x_i+\epsilon)\cap F = \{x_i\}$, then by Corollary \ref{cor: cor_sky_2}
$$
E_X(i_{x_i}(D_i),B) \cong E(D_i,B(x_i-\epsilon, x_i+\epsilon)).
$$
Assume first that $\gamma_{i,0}'$ is the identity map (the case where $\gamma_{i,1}'$ is the identity is entirely similar), and consider the short exact sequence
$$
\ses{B(x_i,x_i+\epsilon)}{B(x_i-\epsilon,x_i+\epsilon)}{B(x_i-\epsilon,x_i]}.
$$
Since $B(x_i-\epsilon,x_i] \cong C_0(x_i-\epsilon,x_i]\otimes H_i$, which is a cone, it follows that
\begin{equation*}
\begin{split}
E(D_i,B(x_i-\epsilon,x_i+\epsilon)) & \cong E(D_i,B(x_i,x_i+\epsilon)) \\
& \cong E(D_i,SH_{i+1}) = 0
\end{split}
\end{equation*}
since $D_i, H_{i+1} \in \mathcal{F}$
\end{proof}

Recall that if $B \in \mathcal{E}(\mathcal{F})$, then by Definition~\ref{def:elementary-1}
for each $i \in \{1,2,\ldots, n\}$, either $\gamma'_{i,0}$ or $\gamma'_{i,1}$ is the identity map.
The corresponding index is denoted by $j(i)$ and $j'(i)=1-j(i)$. In particular this means that $H'_{j(i)}=D'_i$ and $\gamma'_{i,j(i)}=\mathrm{id}$.

\begin{thm} \label{thm: ex}
If $A,B \in \mathcal{E}(\mathcal{F})$ with $\mathcal{F}$ as in Definition~\ref{defn:classF}, then
$$
E_X(A,B) = \left\lbrace (\beta_i) \in \bigoplus_{i=1}^{n+1} E(H_i,H_i') : \beta_{i+j'(i)}\circ\class{\gamma_{i,j'(i)}} = \class{\gamma_{i,j'(i)}'}\circ \beta_{i+j(i)}\circ \class{\gamma_{i,j(i)}} \right\rbrace
$$
\end{thm}
\begin{proof}
By Lemma \ref{lem: af_b} and the exact sequence \eqref{eqn: les_eab}, we see that
$$
E_X(A,B) \cong \ker(\delta : E_X(A(U),B) \to E_X^1(A(F),B)).
$$
But $\ker(\delta)$ was computed in Theorem~\ref{thm: delta}. We deduce that
\begin{equation}\label{eq:calculation}
\begin{split}
E_X(A,B) = \begin{Bmatrix}
				(\beta_i) \in \bigoplus_{i=1}^{n+1} E(H_i,H_i') : & \exists (\alpha_i) \in \bigoplus_{i=1}^n E(D_i,D_i') \text{ s.t.} \\
				& \beta_i\circ \class{\gamma_{i,0}} = \class{\gamma_{i,0}'}\circ \alpha_i \text{ and } \\
				& \beta_{i+1} \circ \class{\gamma_{i,1}} = \class{\gamma_{i,1}'} \circ \alpha_i
			\end{Bmatrix}.
			\end{split}
\end{equation}
The various maps in this description of $E_X(A,B)$ are illustrated in the diagram below.
\[
  \begin{xy}
   \xymatrix{
H_i\ar[ddd]_{\beta_i} &{} & H_{i+1}\ar[ddd]_{\beta_{i+1}}  &{} & H_{i+2}\ar[ddd]^{\beta_{i+2}} \\
{} & D_i\ar[lu]_{\gamma_{i,0}}\ar[ru]^{\gamma_{i,1}} \ar@{-->}[d]^{\alpha_i}&{}&{D_{i+1}}\ar[lu]_{\gamma_{i+1,0}}\ar[ru]^{\gamma_{i+1,1}}\ar@{-->}[d]^{\alpha_{i+1}} &{}\\
{} & D'_i\ar[ld]^{\gamma'_{i,0}}\ar[rd]_{\gamma'_{i,1}} &{}&{D'_{i+1}}\ar[ld]^{\gamma'_{i+1,0}}\ar[rd]_{\gamma'_{i+1,1}} &{}\\
H'_i &{} & H'_{i+1}  &{} & H'_{i+2} \\
   }
  \end{xy}
\]
Let us note that the equations
\begin{equation}\label{eq:calculation++}
\beta_i\circ \class{\gamma_{i,0}} = \class{\gamma_{i,0}'}\circ \alpha_i, \qquad
\beta_{i+1} \circ \class{\gamma_{i,1}} = \class{\gamma_{i,1}'} \circ \alpha_i
\end{equation}
determine $\alpha_i$ uniquely since either $\class{\gamma'_{i,0}}=\mathrm{id}$ or $\class{\gamma'_{i,1}}=\mathrm{id}$.
Indeed, using the notation introduced above, we deduce from \eqref{eq:calculation++}
that $\alpha_i=\beta_{i+j(i)}\circ \class{\gamma_{i,j(i)}}$. Then we substitute this expression in the equation
$\beta_{i+j'(i)}\circ \class{\gamma_{i,j'(i)}}=\class{\gamma_{i,j'(i)}}\circ \alpha_i$ to obtain that
\begin{equation}\label{eq:calculation+}
\beta_{i+j'(i)}\circ\class{\gamma_{i,j'(i)}} = \class{\gamma_{i,j'(i)}'}\circ \beta_{i+j(i)}\circ \class{\gamma_{i,j(i)}}.
\end{equation}
 Conversely, if  \eqref{eq:calculation+} is satisfied for all $i$, then $\alpha_i:=\beta_{i+j(i)}\circ \class{\gamma_{i,j(i)}}$ will satisfy both equations from \eqref{eq:calculation++}.
\end{proof}
\begin{cor}\label{cor:pullback}
Let $Y, Z$ be two closed sub-intervals of $[0,1]$ such that their endpoints are not in $F$. Then $E_{Y\cup Z}(A(Y\cup Z), B(Y\cup Z))$ is the pullback of the following diagram
$$
\begin{xy}
\xymatrix{
E_{Y\cup Z}(A(Y\cup Z), B(Y\cup Z))\ar@{->}[r]\ar@{->}[d] & E_Y(A(Y),B(Y))\ar@{->}[d] \\
E_Z(A(Z),B(Z))\ar@{->}[r] & E_{Y\cap Z}(A(Y\cap Z), B(Y\cap Z))
}
\end{xy}
$$
\end{cor}
\begin{proof} If $I=[a,b]$ is  a closed sub-interval of $[0,1]$ such that its endpoints are not in $F$,
let $i^0_I,i^1_I$ be uniquely defined  by the requirement that $a\in U_{i^0_I}$ and $b\in U_{i^1_I}$.
Let $Y$ and $Z$ be as in the statement. If $Y\cap Z=\emptyset$ the result follows from Theorem~\ref{thm: ex}. Thus we may assume that $Y\cap Z\neq \emptyset$ and moreover that
$i^0_Y\leq i^0_Z \leq i^1_Y \leq i^1_Z$. In this case $i^0_{Y\cap Z}=i^0_Z$,
$i^1_{Y\cap Z}=i^1_Y$ and $i^0_{Y\cup Z}=i^0_Y$,
$i^1_{Y\cup Z}=i^1_Z$.
The statement follows now immediately, since by Theorem~\ref{thm: ex} we have that for each sub-interval $I$ as above
$$
E_I(A,B) = \left\lbrace (\beta_i)_{i^0_I\leq i\leq i^0_I} :  (\beta_i) \ \text{satisfy}\  \eqref{eq:calculation+}\right\rbrace.\qedhere
$$
\end{proof}


\section{Morphisms of the K-theory sheaf}\label{sec: k_theory}

In this section, we apply Theorem \ref{thm: ex} to compute the group $E_{[0,1]}(A,B)$ using K-theory. More precisely, we show that if $A$ and $B$ are elementary $C[0,1]$-algebras whose fibers satisfy the UCT and have $K_0$-groups that are free of finite rank and zero $K_1$-groups, then there is a natural isomorphism
$$
E_{[0,1]}(A,B) \cong \hombbk{A}{B}
$$
where $\mathbb{K}_0(\cdot)$ denotes the K-theory pre-sheaf, an invariant for $C[0,1]$-algebras introduced  in \cite{mdd_elliott}. As an application, we give a partially new proof of the main classification result of \cite{mdd_elliott} which does not require two technical results, Theorem 2.6 and Theorem 8.1, from \cite{mdd_elliott} and instead it uses results of Kirchberg \cite{kirchberg}. \\

We recall the following definition from \cite[$\mathsection 4$]{mdd_elliott}.
\begin{defn}
Let $X$ denote the unit interval and let $A$ be a $C[0,1]$-algebra. Let $\mathcal{I}$ denote the set of all closed subintervals of $X$ with positive length. To each $I \in \mathcal{I}$, associate the group $K_0(A(I))$, and to each pair $I, J \in \mathcal{I}$ such that $J \subset I$, associate the map $$r_J^I = K_0(\pi_J^I): K_0(A(I)) \to K_0(A(J))$$ where $\pi_J^I : A(I) \to A(J)$ is the natural projection. \\

This data gives a pre-sheaf on $\mathcal{I}$ which is denoted by $\mathbb{K}_0(A)$. \\

A morphism of pre-sheaves $\f : \mathbb{K}_0(A) \to \mathbb{K}_0(B)$ consists of a family of maps $\f_I : K_0(A(I)) \to K_0(B(I)))$ such that the following diagram commutes
\begin{equation} \label{eqn: k_map}
\begin{CD}
K_0(A(I)) @> \f_I >> K_0(B(I)) \\
@V r^I_J VV	@VV r^I_J V \\
K_0(A(J)) @> \f_J >> K_0(B(J))
\end{CD}
\end{equation}
  The set of all such morphisms, denoted  $\hombbk{A}{B}$, has an abelian group structure.
\end{defn}

Note that, by \cite[Proposition 2.31]{mdd_meyer}, for each $I \in \mathcal{I}$, there is a natural restriction map
$$
E_X(A,B) \to E_I(A(I),B(I)) \to E(A(I),B(I))
$$
Multiplying $K_0(A(I)) = E(\mathbb{C},A(I))$ with $E(A(I),B(I))$ gives a map
$$
E_X(A,B) \to \homk{A(I)}{B(I)}
$$
Furthermore, if $J \subset I$, then the naturality of the restriction map ensures that the diagram \eqref{eqn: k_map} commutes. Hence, we have a natural map
$$
\Gamma_{A,B} : E_X(A,B) \to \hombbk{A}{B}.
$$

\begin{defn}\label{fdef:F0}
We now introduce a class of algebras for which this map is an isomorphism. Let $\mathcal{F}_0$ be the class of separable  C*-algebras $D$ satisfying the UCT and such that $K_0(D)$ is free of finite rank, and $K_1(D) = 0$.
We define $\mathcal{E}(\mathcal{F}_0)$ to be the class of all elementary $C[0,1]$-algebras whose fibers lie in $\mathcal{F}_0$.
\end{defn}
 \begin{remark}Let us note that the UCT implies that if $D,H \in \mathcal{F}_0$,
then $E^1(D,H) = 0$ and hence that $\mathcal{F}_0\subset \mathcal{F}$.
Thus the results from the previous section apply to members of $\mathcal{E}(\mathcal{F}_0)$.

 Furthermore, for any $H,H' \in \mathcal{F}_0$, the UCT gives us an isomorphism
\begin{equation} \label{eqn: e_k_iso}
E(H,H') \cong KK(H,H') \xrightarrow{} \homk{H}{H'}.
\end{equation}
\end{remark}
Our goal is to show that the map $\Gamma_{A,B}$ is an isomorphism if $A, B \in \mathcal{E}(\mathcal{F}_0)$. In order to do this, we begin by choosing a subset of closed intervals which, roughly speaking, will allow us to capture the K-theory pre-sheaf from a finite amount of data: For each $i \in \{1,2,\ldots, n\}$, choose closed subintervals $V_{i,0}\subset (x_{i-1}, x_i]$ and $V_{i,1}\subset [x_i, x_{i+1})$ both containing $x_i$ and such that
 $V_i = V_{i-1,1}\cap V_{i,0} $ is a nondegenerate interval.
 Using the notation from Theorem~\ref{thm: ex}, we consider the group
$$
G(A,B) = \{(\f_i) \in \homk{A(V_i)}{B(V_i)} : \f_{i+j'(i)}\circ [\gamma_{i,j'(i)}] = [\gamma_{i,j'(i)}']\circ \f_{i+j(i)}\circ [\gamma_{i,j(i)}] \}
$$
Here $[\gamma_{i,j}]$ stands for $K_0(\gamma_{i,j}):K_0(D_i)\to K_0(H_{i+j})$.
\begin{lemma} \label{rem: theta}
There is an isomorphism of groups $\theta : E_X(A,B) \to G(A,B)$
\end{lemma}
\begin{proof} Since each $V_i$ is a closed interval and $A(V_i)=C(V_i, H_i)$, $B(V_i)=C(V_i,H_i')$
we can identify  $\homk{A(V_i)}{B(V_i)}$ with $\homk{H_i}{H'_i}$.
The result  follows now  from Theorem \ref{thm: ex} and the isomorphism \eqref{eqn: e_k_iso}.
The map $\theta$ is induced by the functor that takes an E-theory element to the morphism that it induces
on K-theory groups.
\end{proof}

We now construct a map $\nu : \hombbk{A}{B} \to G(A,B)$ such that
$\nu\circ \Gamma_{A,B} = \theta$.
\begin{lemma} \label{lem: nu_injective}
The map $\nu : \hombbk{A}{B} \to G(A,B) $ given by $\f \mapsto (\f_{V_i})_{i=1}^{n+1}$ is well-defined and injective.
\end{lemma}
\begin{proof}
For any closed interval $I=[a,b]\subset (x_{i-1}, x_{i+1})$ with $x_i\in I,$ we use the extension
$$
\ses{C_0[a, x_i)\otimes H_{i}\oplus C_0(x_i, b]\otimes H_{i+1}}{A(I)}{D_i},
$$
to see that $K_0(A(I))\cong K_0(D_i)$. A similar argument for $B$ shows that $K_0(B(I))\cong K_0(D'_i)$.
In particular $K_0(A(V_{i,0}))\cong K_0(D_i) \cong K_0(A(V_{i,1}))$ and $K_0(B(V_{i,0}))\cong K_0(D'_i) \cong K_0(B(V_{i,1})).$

It follows that if $\f \in \hombbk{A}{B}$, then we can identify the two maps $\f_{V_{i,0}} =  \f_{V_{i,1}} : K_0(D_i) \to K_0(D_i')$.
On the other hand consider the inclusion $V_{i+1} \subset V_{i,1}$, and note that $K_0(A(V_{i+1})) \cong K_0(H_{i+1})$.
We now see  that the restriction map
$$
r^{V_{i,1}}_{V_{i+1}}:K_0(A(V_{i,1})) \to K_0(A(V_{i+1}))
$$
is given by $[\gamma_{i,1}]$. A similar property holds for $B$.
Since any $\f \in \hombbk{A}{B}$ is compatible with the restriction maps, the following diagram is commutative.
$$
\begin{CD}
K_0(A(V_{i,1})) @> [\gamma_{i,1}] >> K_0(A(V_{i+1})) \\
@V \f_{V_{i,1}} VV @V V\f_{V_{i+1}} V \\
K_0(B(V_{i,1})) @> [\gamma_{i,1}'] >> K_0(B(V_{i+1}))
\end{CD}
$$
Thus,
\begin{equation} \label{eqn: alpha_i_1}
\f_{V_{i+1}}\circ [\gamma_{i,1}] = [\gamma_{i,1}']\circ \f_{V_{i,1}}.
\end{equation}
Applying the same argument to the inclusion $V_i \subset V_{i,0}$, we get
\begin{equation} \label{eqn: alpha_i_0}
\f_{V_i}\circ [\gamma_{i,0}] = [\gamma_{i,0}']\circ \f_{V_{i,0}}.
\end{equation}
We saw that $\f_{V_{i,0}} =  \f_{V_{i,1}} : K_0(D_i) \to K_0(D_i')$.
Since $\gamma'_{i,j(i)}=\mathrm{id}$, it follows from \eqref{eqn: alpha_i_1} and \eqref{eqn: alpha_i_0} that
$$
\f_{V_{i+j'(i)}}\circ [\gamma_{i,j'(i)}] = [\gamma_{i,j'(i)}']\circ \f_{V_{i+j(i)}}\circ [\gamma_{i,j(i)}].
$$
This shows that $\nu$ is well-defined.

Now to prove injectivity, suppose that $\f_{V_i} = 0$ for all $i$. We need to show that $\f_I = 0$ for any non-degenerate interval $I \subset [0,1]$. \\
Suppose first that $I$ contains at most one point of $F$. If $I\cap F = \emptyset$, then  $\f_I = \f_{V_i}$ for some $i$, and there is nothing to prove. So suppose that $x_i \in I$ and that no other point of $F$ is in $I$. In that case, $K_0(A(I)) \cong K_0(A(x_i)) \cong K_0(D_i)$
and $\f_I$ can be identified with both $\f_{V_{i,0}}$ and $\f_{V_{i,1}}$,  as seen earlier in the proof. Hence, one of the equations \eqref{eqn: alpha_i_1} or \eqref{eqn: alpha_i_0} (depending on which $\gamma_{i,j}'$ is the identity map) would ensure that $\f_I = 0$. \\

Now consider any nondegenerate interval $I\subset [0,1]$ with $|I\cap F| \geq 2$, and write $I = I_1\cup I_2$, where $I_1$ and $I_2$ are two closed intervals such that $I_1\cap I_2 = \{x\}$ and $x \notin F$, and $I_1$ contains exactly one point of $F$. Then, $A(I)$ can be described as a pull-back
$$
\begin{CD}
A(I) @>>> A(I_1) \\
@VVV @VVV \\
A(I_2) @>>> A(x)
\end{CD}
$$
Applying the Mayer-Vietoris sequence in K-theory, and using the fact that $K_1(A(x)) = 0$, we see that
$$
\begin{CD}
0 @>>> K_0(A(I)) @>>> K_0(A(I_1))\oplus K_0(A(I_2)) @>>> K_0(A(I_1\cap I_2))\\
@. @V \f_{I} VV @VV \f_{I_1}\oplus \f_{I_2} V @VVV\\
0 @>>> K_0(B(I)) @>>> K_0(B(I_1))\oplus K_0(B(I_2)) @>>>K_0(A(I_1\cap I_2))
\end{CD}
$$
Hence, it follows that if $\f_{I_1} = \f_{I_2} = 0$, then $\f_I = 0$. We know from the first part that $\f_{I_1} = 0$, so it suffices to prove that $\f_{I_2} = 0$. We then break up $I_2$ as we did with $I$ before and repeat the same process until we reach an $I_k$ such that $I_k$ contains at most one point of $F$, in which case $\f_{I_k} = 0$ and we can stop the inductive process. This proves the injectivity of $\nu$.
\end{proof}

\begin{thm} \label{thm: uct}
If $A,B \in \mathcal{E}(\mathcal{F}_0)$ (see Def.~\ref{fdef:F0}), then $\Gamma_{A,B} : E_X(A,B) \to \hombbk{A}{B}$ is an isomorphism.
\end{thm}
\begin{proof}
The maps
$
\theta : E_X(A,B) \to G(A,B)$ and $\nu : \hombbk{A}{B} \to G(A,B)
$
satisfy
$\nu\circ \Gamma_{A,B} = \theta$.
By Lemma \ref{rem: theta}, $\theta$ is bijective, and hence $\Gamma_{A,B}$ is injective. By Lemma \ref{lem: nu_injective}, $\nu$ is injective, and hence $\Gamma_{A,B}$ is surjective as well.
\end{proof}

Let $A$ be a separable continuous field over $[0,1]$ whose fibers have vanishing $K_1$ groups. By \cite[Proposition 4.1]{mdd_elliott}, $\mathbb{K}_0(A)$ is a sheaf.
We shall use Theorem~\ref{thm: uct} to give a streamlined proof of the main result of \cite{mdd_elliott}, see Theorem~\ref{thm:DE}.
\emph{For the remainder of this section we make the blanket assumption that all the continuous fields that we consider
are separable and their fibers are stable Kirchberg C*-algebras  with  vanishing $K_1$-groups.}\\

Our definition of elementary $C[0,1]$ algebras given in Def.~\ref{def:elementary-1} is a bit more general that the definition of elementary algebras in the sense of \cite{mdd_elliott}. To make the distinction we will call the latter special elementary. Suppose that $A$ is a special elementary continuous field of Kirchberg algebras. This means that $A$ is defined as the pullback of a certain diagram $\mathcal{D}$. Here is the description of $\mathcal{D}$ where we adapt the notation from \cite{mdd_elliott} to the present setting. \\

Consider a partition $0 = x_0 < x_1 < \ldots < x_{n+1} = 1$, where $n=2m$. Let $A$ be as in Def.~\ref{def:elementary-1}, but we require that
$
D_{2i-1} = H_{2i} = D_{2i}$ and $\gamma_{2i-1,1} = \text{id} = \gamma_{2i,0}.
$
Set
$Y = [x_0,x_1]\cup [x_2,x_3]\cup \ldots \cup [x_{2m},x_{2m+1}]$,
$Z = [x_1,x_2]\cup [x_3,x_4]\cup \ldots \cup [x_{2m-1}, x_{2m}]$
and
$F = \{x_1, x_2, \ldots, x_{2m}\} = Y\cap Z.$
Define
$$
H = \bigoplus_{i=0}^m C[x_{2i},x_{2i+1}]\otimes H_{2i+1}\quad \text{ and }\quad D = \bigoplus_{i=1}^m C[x_{2i-1},x_{2i}]\otimes D_{2i}
$$
whence
$$
H(F) = H_1\oplus \left( \bigoplus_{i=2}^m (H_{2i-1}\oplus H_{2i-1}) \right ) \oplus H_{2m+1}\quad \text{and}\quad  D(F)=\bigoplus_{i=1}^{2m }D_i.
$$
Consider the diagram $\mathcal{D}$ given by
$$
\begin{CD}
H @>{\pi}>> H(F) @< \gamma << D
\end{CD}
$$
where $\pi$ is the restriction map, and $\gamma$ is the composition of the map $\gamma':D(F)\to H(F)$, with components $\gamma_{2i-1,0}:D_{2i-1}\to H_{2i-1}$, $\gamma_{2i,1}:D_{2i}\to H_{2i+1}$, with the restriction map $D\to D(F)$. $A$ is isomorphic to the pullback of the diagram $\mathcal{D}$ and we have an induced commutative diagram
$$
\begin{xy}
\xymatrix{
A(Y)\ar[r]^{\pi}\ar[d] & A(F)\ar[d]_{\gamma'}& A(Z)\ar[l]_{\pi}\ar@{=}[d] \\
H\ar[r] & H(F) &D\ar[l]_{\gamma}
}
\end{xy}
$$
Given $\mathcal{D}$ as above, and $B$  a continuous field over $[0,1]$, we denote by $\mathcal{D}B$ the diagram
$$
\begin{CD}
B(Y) @>{\pi}>> B(F) @< {\pi} << B(Z).
\end{CD}
$$
Recall from \cite{{mdd_elliott}} that a fibered morphism  $\varphi \in \text{Hom}_{\mathcal{D}}(A, B)$
is given a commutative diagram
$$
\begin{CD}
H @>>> H(F) @<<< D \\
@V \varphi_Y VV @V \varphi_F VV @VV \varphi_Z V \\
B(Y)@>>> B(F) @<<< B(Z).
\end{CD}
$$
where the vertical arrows are injective monomorphisms of continuous fields.
The combination of the two larger diagram above gives a morphism of diagrams
$\mathcal{D}A \to \mathcal{D}B$ which induces a fiberwise injective morphism of continuous
fields $\widehat{\varphi}:A \to B$. A morphism of fields induced by a fibered morphism is called \emph{elementary}
 \cite[p.806]{mdd_elliott}.
As in \cite{mdd_elliott}, denote by $K_0(\mathcal{D})$, the diagram
$$
\begin{CD}
K_0(H) @>{\pi_{\ast}}> >K_0(H(F)) @<{\gamma_{\ast}} <<K_0(D).
\end{CD}
$$
$\text{\textnormal{Hom}}(K_0(\mathcal{D}), K_0(\mathcal{D}B))$ consists of all morphisms of diagrams of groups
$$
\begin{CD}
K_0(H) @>>> K_0(H(F)) @<<< K_0(D) \\
@V \alpha_Y VV @V \alpha_F VV @VV \alpha_Z V \\
K_0(B(Y)) @>>> K_0(B(F)) @<<< K_0(B(Z))
\end{CD}
$$
that  preserve the direct sum decomposition of the K-theory groups induced by the underlying partition of $[0,1]$.
It is a K-theory counterpart of $\text{Hom}_{\mathcal{D}}(A, B)$.
In \cite[Prop.5.1]{mdd_elliott}, it is shown that each $\alpha \in \text{Hom}(K_0(\mathcal{D}), K_0(\mathcal{D}B))$ induces a morphism of sheaves
$\widehat{\alpha} \in \hombbk{A}{B}.
$
Let us also note that a morphism $\beta\in \hombbk{B}{B'}$
induces by restriction  a morphism $\mathcal{D}\beta\in \text{\textnormal{Hom}}(K_0(\mathcal{D}B), K_0(\mathcal{D}B'))$.
To simplify notation, we will often write $\beta$ in place of  $\mathcal{D}\beta$.

\begin{lemma}{\cite[Prop.5.1]{mdd_elliott}} \label{lemma: elementary_diagram}
\begin{itemize}
\item[(i)] Suppose that $K_0(H_i)$ and $K_0(D_i)$ are finitely generated.
If $B = \varinjlim B_n$ is an inductive limit of continuous fields $B_n$ over $[0,1]$, then
 \[\text{Hom}(K_0(\mathcal{D}), K_0(\mathcal{D}B)) \cong \varinjlim\text{Hom}(K_0(\mathcal{D}), K_0(\mathcal{D}B_n))\]
\item[(ii)] If $\alpha \in \text{\textnormal{Hom}}(K_0(\mathcal{D}), K_0(\mathcal{D}B))$ and  $\beta \in \hombbk{B}{B'}$, then
    $((\mathcal{D}\beta)\circ \alpha)\,\widehat{}=\beta\circ\widehat{\alpha}$.
\end{itemize}
\end{lemma}
\begin{proof}
(i)  Since $\mathcal{D}$ is a finite diagram of finitely generated groups, the result follows using the continuity of the $K_0$-functor.
 (ii) is proved in \cite[Prop.5.1]{mdd_elliott}.
\end{proof}

We are now ready to reprove the classification theorem of \cite{mdd_elliott}:
\begin{thm}\label{thm:DE}
Let $A,B$ be separable continuous fields over $[0,1]$ whose fibers are stable Kirchberg C*-algebras satisfying  the UCT, with torsion free $K_0$-groups  and vanishing $K_1$-groups. If $\alpha \in \hombbk{A}{B}$ is an isomorphism of sheaves, then there is an isomorphism of continuous fields $\phi : A \to B$ such that $\phi_{\ast} = \alpha$
\end{thm}
\begin{proof}
Recall that the inductive limit decomposition from Theorem~\ref{thm:structure} comes with more structure
and properties that we now review. Specifically, in the inductive system
$$
A_1 \xrightarrow{\widehat{\varphi}_{1,2}} A_2 \xrightarrow{\widehat{\varphi}_{2,3}} \ldots \to A_n \xrightarrow{\widehat{\varphi}_{n,n+1}} A_{n+1} \to \ldots
$$
with $A\cong \varinjlim A_k$, all connecting morphisms are elementary in the sense of \cite[p.806]{mdd_elliott}. In other words, each $A_k$ is the pull-back of a diagram $\mathcal{D}_k$, and each $\widehat{\varphi}_{k,k+1} \in \text{Hom}(A_k,A_{k+1})$ is induced by a fibered morphism  $\varphi_{k,k+1}\in \text{Hom}_{\mathcal{D}_k}(A_k,A_{k+1})$. Moreover ${\varphi}_{k,\infty} \in \text{Hom}_{\mathcal{D}_k}(A_k,A)$ and
$ \widehat{\varphi}_{k+1,\infty}\circ {\varphi}_{k,k+1}={\varphi}_{k,\infty}$.
The fibers of $A_k$ are stable Kirchberg algebras whose $K_0$-groups are free of finite rank  and their $K_1$-groups vanish.
Similarly, let $(B_n)$ be a sequence approximating $B$, where $B_n$ is the pull-back of the diagram $\mathcal{D}_n'$, and ${\psi}_{n,n+1}\in \text{Hom}_{\mathcal{D}_n'}(B_n, B_{n+1}), {\psi}_{n,\infty} \in \text{Hom}_{\mathcal{D}_n'}(B_n,B)$ are the corresponding maps. \\

 By Lemma \ref{lemma: elementary_diagram}(i), for $\alpha_1: = \alpha \circ (\varphi_{1,\infty})_{\ast}\in \text{Hom}(K_0(\mathcal{D}_1), K_0(\mathcal{D}_1 B))$,  there is $m_1 \in \mathbb{N}$ and $\mu_1 \in \text{Hom}(K_0(\mathcal{D}_1), K_0(\mathcal{D}_1B_{m_1}))$ such that $\alpha_1 = (\widehat{\psi}_{m_1,\infty})_{\ast}\circ \mu_1$. Letting $n_1=1$ we have $\widehat{\mu_1} \in \hombbk{A_{n_1}}{B_{m_1}}$.
Similary, for $\beta_1 := \alpha^{-1}\circ (\psi_{m_1,\infty})_{\ast} \in \text{Hom}(K_0(\mathcal{D}_{m_1}'),K_0(\mathcal{D}_{m_1}'A))$, there is $n_2 >n_1$  and $\eta_1 \in \text{Hom}(K_0(\mathcal{D}_{m_1}'), K_0(\mathcal{D}_{m_1}A_{n_2}))$ such that $\beta_1 = (\widehat{\varphi}_{n_2,\infty})_{\ast}\circ \eta_1$. This gives $\widehat{\eta}_1 \in \hombbk{B_{n_1}}{A_{n_2}}$.
Combining the equations $\alpha \circ (\varphi_{n_1,\infty})_{\ast}=(\widehat{\psi}_{m_1,\infty})_{\ast}\circ \mu_1$
and $ \alpha^{-1}\circ (\psi_{m_1,\infty})_{\ast}=(\widehat{\varphi}_{n_2,\infty})_{\ast}\circ \eta_1$, we use
 Lemma \ref{lemma: elementary_diagram}(ii) to deduce that
 $(\widehat{\varphi}_{n_2,\infty})_\ast \circ \widehat{\eta}_1\circ \mu_1=(\widehat{\varphi}_{n_1,\infty})_\ast=(\widehat{\varphi}_{n_2,\infty})_\ast \circ ({\varphi}_{n_1,n_2})_\ast$. By  Lemma \ref{lemma: elementary_diagram}(i)
 we conclude that after increasing $n_2$, if necessary, we can arrange that
$\widehat{\eta}_1\circ \mu_1=  ({\varphi}_{n_1,n_2})_\ast$ and hence $\widehat{\eta}_1\circ \widehat{\mu}_1=  (\widehat{\varphi}_{n_1,n_2})_\ast.$
By induction, we construct a commutative diagram
\[
  \begin{xy}
   \xymatrix{
\mathbb{K}_0(A_{n_1}) \ar[dr]^{\widehat{\mu}_1} \ar[rr]^{(\widehat{\varphi}_{n_1,n_2})_{\ast}} &{} & \mathbb{K}_0(A_{n_2}) \ar[dr]\ar[rr] &{} & \mathbb{K}_0(A_{n_3}) \ar[r] & \ldots \ar[r] & \mathbb{K}_0(A) \ar@<1.ex>[d]^{\alpha} \\
& \mathbb{K}_0(B_{m_1}) \ar[ru]^{\widehat{\eta}_1}\ar[rr]^{(\widehat{\psi}_{m_1,m_2})_{\ast}}&{} & \mathbb{K}_0(B_{m_2}) \ar[ru]\ar[rr] &{} & \ldots \ar[r] & \mathbb{K}_0(B)\ar@<1.ex>[u]^{\alpha^{-1}}
   }
  \end{xy}
\]

By Theorem \ref{thm: uct}, we replace the diagonal arrows by $E_X$-theory elements and hence by  $KK_X$-elements, since all involved continuous fields are nuclear \cite{park_trout}.
 By Kirchberg's results \cite{kirchberg}, we can further replace these $KK_X$-elements by morphisms of fields which are fiberwise injective and moreover, each triangle is commutative up to asymptotic unitary equivalence. This yields an isomorphism $\phi: A\to  B$ with $\phi_\ast=\alpha$,  by applying Elliott's intertwining argument \cite[Sec.~2.3]{rordam_stormer}.
\end{proof}
\bibliographystyle{plain}
\bibliography{mybib}

\end{document}